\documentclass[notnumberedtheorems,withtitlethanks]{actacybpress}
\usepackage{graphicx, algorithm, algorithmic}

\begin{document}

% Overfull boxes can be avoided by using the \sloppy command.
% Use it as a last chance only since it produces underfull boxes.
% \fussy resets to the normal setting.
%\sloppy
%...
%\fussy

%% Uncomment commands below once the required information is known
%% Used by the technical editor
%\acta{vol}{num}{year}{paper_id}
%\setcounter{page}{163}
%\received{...} 

%\title{SCAN 2020: Affine Iterations and Wrapping Effect: Various Approaches}
% for blind submission
%\author{~ \thanks{~ \email{~} \orcid{~}}}
% for HAL 
\title{Affine Iterations and Wrapping Effect: Various Approaches}
\headingtitle{Affine Iterations and Wrapping Effect}

\author{Nathalie Revol
\thanks{INRIA - LIP, ENS de Lyon, France,
\email{Nathalie.Revol@inria.fr},
\orcid{https://orcid.org/0000-0002-2503-2274}}}

% Please provide ORCIDs for all authors!

\maketitle

\begin{abstract}
Affine iterations of the form \(x_{n+1}=Ax_n+b\) converge, using real arithmetic, if the spectral radius of the matrix \(A\) is less than 1. However, substituting interval arithmetic to real arithmetic may lead to divergence of these iterations, in particular if the spectral radius of the absolute value of \(A\) is greater than 1. 
%We will recall some results due to Mayer and co-authors regarding the divergence of the iterates.
We will review  different approaches to limit the overestimation of the iterates, when the components of the initial vector \(x_0\) and \(b\) are intervals. 
%and we will propose a new approach based on the SVD decomposition of \(A\).
We will compare, both theoretically and experimentally, the widths of the iterates computed by these different methods: the naive iteration, methods based on the QR- and SVD-factorization of \(A\), and Lohner's QR-factorization method.
The method  based on the SVD-factorization is computationally less demanding and gives good results when the matrix is poorly scaled, it is superseded either by the naive iteration or by Lohner's method otherwise.

% Put keywords appropriate to your paper here, as shown

% Put your AMS subject classifications into the argument of the following command.
%\AMSsubj{65G20}
\end{abstract}

\section{Introduction} \label{sec:introduction}

The problem we consider is the evaluation of the successive iterates of
\[
\left\{ \begin{array}{lcl} 
x_{n+1} & = & A x_n + b, \\
\multicolumn{3}{l}{x_0 \mbox{ given},}
\end{array} \right.
\]
where \( A \in \mathbb{R}^{d \times d}\), \(x_n\in\mathbb{R}^d\) for every \(n\in\mathbb{N}\) and \(b\in\mathbb{R}^d\).
More specifically, the focus is on the use of interval arithmetic to evaluate these iterates.

In what follows, interval quantities will be denoted in boldface.

\subsection{A Toy Example}
\label{sec:toy-example}
This problem was brought to us through this example of an IIR (Infinite Impulse Response) linear filter in a state-space form:
\[
x_{n} = 1.8*x_{n-1} - 0.9*x_{n-2} + 4.7.10^{-2}*(u_{n-2}+u_{n-1}+u_{n})
\]
for \( x_0 = 0\) and \(x_1 \in [1, 1.1]\).
We assume that \( u_n \in {\mathbf u} = [ 9.95 \, , \, 10.05]\) for every \(n\).

This iteration can also be written as a linear recurrence in \(\mathbb{R}^2\):
\[ \begin{array}{ll}
& \left( \begin{array}{r} x_{n-1} \\ x_n \end{array} \right)
= A .
\left( \begin{array}{r} x_{n-2} \\ x_{n-1} \end{array} \right) + b_n, \\
\mbox{ where } & A = \left( \begin{array}{rl} 0 & 1 \\ -0.9 & 1.8 \end{array} \right)
\mbox{ and } b_n = \left( \begin{array}{c} 0 \\ 4.7.10^{-2}*(u_{n-2}+u_{n-1}+u_{n}) \end{array} \right) .
\end{array}
 \]

This toy example will be used to illustrate the various approaches mentioned in this paper.
The first  iterates, obtained using floating-point arithmetic, with random values for \(x_1 \in [1, 1.1]\) and each \(u_n \in [ 9.95 \, , \, 10.05] \), are given on the left two columns below. \\
{\small
\[
\begin{array}{cc||cc}
\begin{array}{c|r}
n & x_n \\
\hline
0 & 0 \\
1 & 1.0617 \\
2 & 3.3183 \\
3 & 6.4234 \\
4 & 9.9851 \\
5 & 13.6031 \\
6 & 16.9117 \\
7 & 19.6103 \\
8 & 21.4884 \\
9 & 22.4394 \\
10 & 22.4595 \\
12 & 20.1508 \\
15 & 13.8931 \\
\end{array}
& \hspace*{5mm}
\begin{array}{c|r}
n & x_n \\
\hline
20 & 9.1518 \\
30 & 17.0186 \\
40 & 12.4414 \\
50 & 15.0305 \\
60 & 13.6130 \\
70 & 14.3858 \\
80 & 13.9680 \\
90 & 14.1510 \\
100 & 14.0949 \\
200 & 14.0870 \\
300 & 14.1443 \\
400 & 14.1282 \\
500 & 14.0828 \\
\end{array}
\hspace*{5mm}
& \hspace*{5mm}
\begin{array}{c|r}
n & \mathrm{wid}(\mathbf{x}_n) \\
\hline
0 & 0 \\
1 & 0.1000 \\
2 & 0.1941 \\
3 & 0.4535 \\
4 & 1.0051 \\
5 & 2.2313 \\
6 & 4.9350 \\
7 & 10.905 \\
8 & 24.085 \\
9 & 53.182 \\
10 & 117.42 \\
12 & 572.31 \\
15 & 6158.0 \\
\end{array}
& \hspace*{3mm}
\begin{array}{c|l}
n & \mathrm{wid}(\mathbf{x}_n) \\
\hline
20 & 3.2293. 10^5 \\
30 & 8.8808 . 10^8 \\
40 & 2.4423 . 10^{12} \\
50 & 6.7164 .  10^{15} \\
60 & 1.8470 . 10^{19} \\
70 & 5.0794 . 10^{22} \\
80 & 1.3969 . 10^{26} \\
90 & 3.8415 . 10^{29} \\
100 & 1.0564 . 10^{33} \\
200 & 2.6137 . 10^{67} \\
300 & 6.4663 . 10^{101} \\
400 & 1.5998 . 10^{136} \\
500 & 3.9580 . 10^{170} \\
\end{array}
\end{array}
\]
}

The system stabilizes around \(14\), with variations due to the random values taken by the \(u_n\).
However, the following snippet of Octave code computes the successive iterates using interval arithmetic, using the interval \([1 \, , \,1.1]\) for \({\mathbf x}_1\) and \({\mathbf u} = [9.95 \, , \, 10.05]\) for the \(u_n\), that is, 
%we replace \(4.7.10^{-2}*(u_{n-2}+u_{n-1}+u_{n}) \) by \( {\mathbf b} = 3*4.7.10^{-2} {\mathbf u}\).
we replace \(u_{n-2}+u_{n-1}+u_{n} \) by \(  3* {\mathbf u}\).
\begin{center}
\begin{verbatim}
A=[[0 1];[-0.9 1.8]];
xn=[infsup(0,0);infsup(1,1.1)];
b=4.7e-2 * 3.0*[infsup(0,0);infsup(9.95,10.05)];
n=500; for i=1:n, i , xn=A*xn+b, wid(xn(1)), end;
\end{verbatim}
\end{center}
On the right two columns above are the widths of the successive iterates \({\mathbf x}_n\): the widths of the iterates diverge rapidly to infinity.
%\\

The explanation of this phenomenon is the following: the spectral radius of \(A\) is strictly less than 1: \(\rho(A) \simeq 0.9487 < 1\), and thus the exact (and, for that matter, floating-point) iterations converge.
However, the recurrence satisfied by the widths of the iterates is
\( \mathrm{wid}({\mathbf x}_{n}) = 1.8*\mathrm{wid}({\mathbf x}_{n-1}) + 0.9*\mathrm{wid}({\mathbf x}_{n-2}) + 4.7.10^{-2}*3*\mathrm{wid}({\mathbf u}) \),
which corresponds to the 2-dimensional iteration 
\( w_n = |A|. w_{n-1} + w_b\),
with \(w_n = \mathrm{wid}({\mathbf x}_{n}) \),
\(|A|\) the matrix whose coefficients are the absolute values of the coefficients of \(A\)
and \(w_b = 4.7.10^{-2}*3*\mathrm{wid}({\mathbf u})\).
As the spectral radius of \(|A|\) is larger than 1, indeed \(\rho(|A|) \simeq 2.208 >1\), the iterations diverge.
%\\

This phenomenon is a special case of the so-called {\it wrapping effect}. Its ubiquity in interval computations has been put in evidence by Lohner in \cite{Lohner2001}.

\subsection{The Wrapping Effect}

The wrapping effect is ubiquitous, as defined and developed in \cite{Lohner2001}.
It can be described as the overestimation due to the enclosure of the sought set in a set of a given. simple structure.
In our case, this simple structure corresponds to multidimensional intervals or boxes, that is, parallelepipeds with sides parallel to the axes of the coordinate system.
When the computation is iterative, and when each iteration produces such an overestimating set that is used as the starting point  of the next iteration, the size of the computed set may grow exponentially in the number of iterations, even when the exact solution set remains bounded and small.

Lohner also put in evidence that the affine iteration we study in this paper, namely \( x_{n+1} = Ax_n + b\), or more generally 
\( x_{n+1} = A_n x_n + b_n\) with \( x_{n+1}\), \(x_n\) and \(b_n\) vectors in \(\mathbb{R}^d\) and \(A_n \in \mathbb{R}^{d \times d}\) for every \(n \in \mathbb{N}\), is archetypal. It occurs in many algorithms, and the examples cited in \cite{Lohner2001} include 
\begin{itemize}
\item matrix-vector iterations   as the ones studied in this paper; 
\item discrete dynamical systems: \( {\mathbf x}_{n+1} = f({\mathbf x}_n)\), \({\mathbf x}_0\) given and \(f\) sufficiently smooth;
\item continuous dynamical systems (ODEs): \( x'(t)= g(t,x(t))\), \(x(0)=x_0\),
which is studied through a numerical one step method (or more) of the kind \( {\mathbf x}_{n+1} = {\mathbf x}_n + h \Phi({\mathbf x}_n,t_n) + {\mathbf z}_{n+1}\);
\item difference equations: \( a_0 {\mathbf z}_n + a_1 {\mathbf z}_{n+1} + \ldots + a_m {\mathbf z}_{n+m} = b_n\) with \({\mathbf z}_1, \ldots {\mathbf z}_m\) given;
\item linear systems with (banded) triangular matrix;
\item automatic differentiation.
\end{itemize}

In this paper, we concentrate on examples similar to the toy example presented above: for every initial value \( x_0 \in {\mathbb R}^d\), the sequence of iterates \( (x_n)_{n \in {\mathbb N}}\) converges to a finite value \(x^* \in {\mathbb R}^d\), since \( \rho(A) <1\);
however, the computations performed using interval arithmetic diverge because their behaviour is dictated by \( \rho (|A|)\) which is larger than \(1\).
We are interested in the iterates computed using interval arithmetic: it is established that these iterates increase in width, however different approaches can be applied to counteract the exponential growth of the width of the iterates.
Several of them, some new  as in  Sections \ref{sec:QR} and \ref{sec:SVD}, and some already well eestablished as in Section~\ref{sec:LohnerQR}, will be tried and compared, in terms of the widths of the results and the computational time.

\section{Theoretical Results}
\label{sec:methods}

\subsection{Problem and Notations}

Let \(A\) be a \(d \times d\) matrix in \(\mathbb{R}^{d \times d}\),
\( {\mathbf x}_0 \in {\mathrm I}\!\mathbb{R} ^d\) be an interval vector (boldface font is used for interval quantities and \( {\mathrm I}\!\mathbb{R}\) stands for the set of real intervals),
\(x_0\) a vector in \(\mathbb{R}^d\) with \(x_0 \in {\mathbf x}_0\), 
\( {\mathbf b} \in {\mathrm I}\!\mathbb{R} ^d\) an interval vector,
and \(b\) a vector in \(\mathbb{R}^d\) and \(b \in {\mathbf  b}\). 
In what follows, \(n\) denotes the number  of iterations.

It is assumed that \( \rho(A)<1\) and \(\rho(|A|)>1\).

A first goal is to determine the set of all fixed-points of the iteration
\[ \left\{ \begin{array}{l}
x_0 \in {\mathbf x}_0, \mbox{ } b \in {\mathbf  b}, \\
x_{n+1} = A x_n + b
\end{array}
\right.
\]
for every \(x_0 \in {\mathbf x}_0\) and every \(b \in {\mathbf b}\).

It is known that \(x_n\) can be written as 
\[ x_n = A^n x_0 + \sum_{i=1}^{n-1} A^i b,\]
thus 
\[ \{ x_n : x_0 \in {\mathbf x}_0 \} \subset A^n {\mathbf x}_0 + \left( \sum_{i=1}^{n-1} A^i \right) {\mathbf b}.\]

However, when the vectors \(x_0\) and \(b\) are replaced in the iterative formula by their interval enclosures \({\mathbf x}_0\) and \({\mathbf b}\), 
one obtains the new interval vector \({\mathbf x}_{n+1}\), which is computed as:
\[ \left\{ \begin{array}{l}
{\mathbf x}_0 \mbox{ and } {\mathbf b} \mbox{ given}, \\
{\mathbf x}_{n+1} = A {\mathbf x}_n + {\mathbf b}.
\end{array}
\right.
\]
Another goal is to determine a tight enclosure for each iterate of this diverging set of intervals.

As mentioned above, the increase in widths of the iterates can be attributed to the use of parallelepipeds with sides parallel to the axes of the coordinate system, and not to the geometry of the transformation.
To cure this problem, changes of coordinates will be applied, using an invertible matrix \(B\), with \(x = B y \Leftrightarrow y = B^{-1}x\) and its interval counterpart \({\mathbf x}=B {\mathbf y}\).
This yields the iteration
\[ \left\{ \begin{array}{lcl}
x_{n+1} & = & By_{n+1} \\
y_{n+1} & = & B^{-1}AB y_n+ B^{-1}b
\end{array}
\right.
\]
and its interval counterpart
\[ \left\{ \begin{array}{lcl}
{\mathbf x}_{n+1} & = & B{\mathbf y}_{n+1} \\
 {\mathbf y}_{n+1} & = & B^{-1}AB  {\mathbf y}_n+ B^{-1}{\mathbf b}.
\end{array}
\right.
\]

In what follows, to establish bounds and their proofs, we assume \(A\) diagonalizable (this will not necessarily be the case for the experiments)
and \(A\) can be diagonalized as \( A = P^{-1} \Lambda P\) where \(\Lambda\) is a diagonal matrix with the eigenvalues \( \lambda_1, \ldots \lambda_d\) of \(A\) on the diagonal and the columns of \(P^{-1}\) are the corresponding eigenvectors.

The iteration considered in this paper corresponds to \( x_n = A^n x_0 + \sum _{i=0}^{n-1} A^i b\). The numerical unstability of computing the matrix power \(A^n\) and applying  it to a vector, is well known:  \(A^n x_0\) tends to be aligned with the eigenvector of \(A\) associated with the largest (in module) eigenvalue, and the information corresponding to the contribution of  the other eigenvectors is lost.
To avoid this well-known problem of the power method, we will consider orthogonal changes of coordinates.
The choice of the orthogonal matrices is related to \(A\), the matrix of the iteration.
\\

We will first consider the QR-factorization of \(A\): \(A = QR\) with \(Q \in {\mathbb R}^{d \times d}\) orthogonal, that is, \( Q Q' = Q' Q = I\) is the identity matrix and \(R \in {\mathbb R}^{d \times d}\) is upper triangular.

The other factorization used in this paper is the SVD-factorization of \(A\): 
\( A = U \Sigma V'\) with \(U\), \(V\) and \(\Sigma\) \( \in {\mathbb R}^{d \times d}\), where \(U\) and \(V\) are orthogonal and \(\Sigma\) is a diagonal matrix with the singular values \( \sigma_1, \ldots \sigma_d\) of \(A\) on the diagonal. We also assume that  \(\sigma_1 \geq \sigma_2 \geq \ldots \geq \sigma_n \geq 0\). This idea has been sketched but not completely developed by Beaumont in \cite{Beaumont2000}.

\subsection{Known results}
Mayer and his co-authors have extensively studied the existence of a fixed-point  for the  iteration studied  in this paper.
In \cite{MayerWarnke2003}, Mayer and Warnke have thoroughly established formulas for the fixed-point  in the case of \( \rho(|A|)<1\): this fixed-point is independent of the starting interval \( {\mathbf x}_0\).
In \cite{ArndtMayer2004}, Arndt and Mayer have  established necessary and sufficient condition on \(A\) for a fixed-point to exist, when \( \rho(|A|)=1\). In this case, the fixed-point is an interval of nonzero width, that is, a non-degenerate interval.
It is well-known that the widths of the iterates diverge when \(\rho(|A|) >1\), and thus that no fixed-point exists  in this case. Our goal is to study the speed of divergence of the iterates in this case.
%From Mayer and Warnke:  no fixed-point of width zero,
%estimation of the width of the \(n\)-th iterate (for comparison with later enclosures).

\subsection{Different Approaches along with Theoretical Bounds}
The main idea is to use an orthogonal change of coordinates which is related to the matrix of the iteration.
As the matrix \(A\) is kept constant for all iterations (and this is not the case in the more general approach of Lohner, see Section~\ref{sec:LohnerQR}), the  change of  coordinates is also kept constant and given by an orthogonal matrix \(B\).
The two orthogonal  matrices considered in what follows are either \(B=Q\) from the QR-factorization of \(A\), or \(B=U\), resp. \(B=V'\), from the SVD-factorization of \(A\).

\subsubsection{Orthogonal Change of Coordinates}
\label{sec:generalB}
Before diving into the  specificities of these changes of coordinates, let us study the general change of coordinates
using an orthogonal matrix \(B\), that is, \( B^{-1} = B'\), with \(x = B y \Leftrightarrow y = B^{-1}x\) and its interval counterpart \({\mathbf x}=B {\mathbf y}\).
%This yields the iteration
%\[ \left\{ \begin{array}{lcl}
%x_{n+1} & = & By_{n+1} \\
%y_{n+1} & = & B^{-1}AB y_n+ B^{-1}b
%\end{array}
%\right.
%\]
The interval iteration  is
\[ \left\{ \begin{array}{lcl}
{\mathbf x}_{n+1} & = & B{\mathbf y}_{n+1} \\
 {\mathbf y}_{n+1} & = & B^{-1}AB  {\mathbf y}_n+ B^{-1}{\mathbf b}.
\end{array}
\right.
\]

Thus the iteration satisfied by the width of the \( {\mathbf y}_n\) is
\[ 
\begin{array}{lcl}
\mathrm{wid}(\mathbf{y}_{n+1}) & = & | B^{-1} A B| \mathrm{wid}(\mathbf{y}_n) + |B^{-1}| \mathrm{wid}(\mathbf{b}) \\
& \leq & |B^{-1}|.|A|.|B| \mathrm{wid}(\mathbf{y}_n) + |B^{-1}| \mathrm{wid}(\mathbf{b})
\end{array}
\]
where the inequalities are to be understood componentwise.
By induction on \(n\),
\[ \mathrm{wid}(\mathbf{y}_n) \leq (|B^{-1}|.|A|.|B|)^n . \mathrm{wid}(\mathbf{y}_0) + \sum_{i=0}^{n-1} (|B^{-1}|.|A|.|B|)^i . |B^{-1}| . \mathrm{wid}(\mathbf{b}).
\]
Taking norms, one gets
\[
\begin{array}{lccl}
\| \mathrm{wid}(\mathbf{y}_n) \| & \leq & & \left( \| \: |B^{-1}| \: \| . \| \: |A|\: \|. \| \: |B|\:\| \right)^n \| . \mathrm{wid}(\mathbf{y}_0)\| \\[2mm]
& & + & \sum_{i=0}^{n-1} \left( \| \: |B^{-1}| \: \|  . \| \: |A|\: \|. \| \: |B|\:\| \right)^i . \| \: |B^{-1}| \: \|. \| \mathrm{wid}(\mathbf{b}) \| \\[4mm]
& \leq & &  \left( \| \: |B^{-1}| \: \| . \| \: |A|\: \|. \| \: |B|\:\| \right)^n \| . \mathrm{wid}(\mathbf{y}_0)\|  \\[2mm]
& & + &  \frac{ \left( \| \: |B^{-1}| \: \| . \| \: |A|\: \|. \| \: |B|\:\| \right)^n -1}{ \| \: |B^{-1}| \: \| . \| \: |A|\: \|. \| \: |B|\:\| -1}  \| \: |B^{-1}| \: \| . \| \mathrm{wid}(\mathbf{b})\|.
\end{array}
\]

\noindent {\bf Remark:} if the considered norm is the matrix norm induced by the vector Euclidean norm, then \( \| \: |B| \: \|_2 = \| B \|_2\) for any matrix \(B\). Similarly, \( \| \: |B| \: \|_{\infty} = \| B \|_{\infty} \leq \sqrt{d}\) for any \(d \times d\) orthogonal matrix \(B\).
In such cases, the bound becomes
\[
\begin{array}{lccl}
\| \mathrm{wid}(\mathbf{y}_n) \| & \leq & & \left( \kappa(B) . \| \: |A|\: \| \right)^n  \| \mathrm{wid}(\mathbf{y}_0)\| \\[2mm]
& & + &  \frac{ \left( \kappa(B) . \| \: |A|\: \|\right)^n -1}{ \kappa(B) . \|     \: |A|\: \| -1}  \| B^{-1} \| . \| \mathrm{wid}(\mathbf{b})\|,
\end{array}
\]
where  \( \kappa(B) \) denotes \( \| B \| . \| B^{-1} \|\), the condition number of \(B\) for the problem of solving a linear system.

Since \( \| B \|_2 = \| B^{-1}\|_2 = \kappa_2(B)=1\) for an orthogonal matrix \(B\),  this bound simplifies even further with the Euclidean norm:
\[
\begin{array}{lcl}
\| \mathrm{wid}(\mathbf{y}_n) \| & \leq &  \|  A \| ^n .  \| \mathrm{wid}(\mathbf{y}_0)\| 
 +   \frac{ \| A \|^n -1}{ \|A \| -1} .  \| \mathrm{wid}(\mathbf{b})\|.
\end{array}
\]
In other words, theoretically there is no difference in the bounds on the widths of the iterates, whether an orthogonal change of coordinates takes place or not.
\\

In what follows, we assume again \(A\) diagonalizable (this will not necessarily be the case for the experiments)
and \(A\) can be diagonalized as \( A = P^{-1} \Lambda P\) where \(\Lambda\) is diagonal.
If we replace \(A\) by \( P^{-1} \Lambda P\) in the iteration, one gets the mathematically equivalent formulation
\[
\begin{array}{lcl}
\mathbf{y}_{n+1} & = & B^{-1} P^{-1} \Lambda P B \mathbf{y}_n + B^{-1} \mathbf{b} \\
& = & (PB)^{-1} \Lambda (PB) \mathbf{y}_n + B^{-1} \mathbf{b},
\end{array}
\]
thus 
\[ \mathrm{wid}(\mathbf{y}_n) = (|(PB)^{-1}|.|\Lambda|.|PB|) . \mathrm{wid}(\mathbf{y}_n) + |B^{-1}|. \mathrm{wid}(\mathbf{b}),
\]
and by induction
\[ \mathrm{wid}(\mathbf{y}_n) = (|(PB)^{-1}|.|\Lambda|.|PB|)^n  . \mathrm{wid}(\mathbf{y}_0) + \sum_{i=0}^{n-1} ( |(PB)^{-1}|.|\Lambda|.|PB|)^i . |B^{-1}|. \mathrm{wid}(\mathbf{b}).
\]
Taking the Euclidean norm of vectors and the induced matrix norm, one gets
\[
\begin{array}{lccl}
\| \mathrm{wid}(\mathbf{y}_n) \|_2 & \leq & & \left( \kappa_2(PB) \| \Lambda\|_2 \right)^n . \| \mathrm{wid}(\mathbf{y}_0)\|_2 \\[2mm]
& & + & \frac{\left( \kappa_2(PB) \| \Lambda\|_2 \right)^n -1}{\kappa_2(PB) \| \Lambda\|_2 -1} . \| \mathrm{wid}(\mathbf{b}) \|_2. \\[4mm]
\end{array}
\]

Let us note that \( \kappa(PB) = \kappa(P) \).
Furthermore, as \( \Lambda\) is diagonal, \( \| \Lambda \| \) is the largest eigenvalue (in module) of \(A\), that is,  \( \| \Lambda \| = \rho(A) <1\). This implies
\[
\begin{array}{lcl}
\| \mathrm{wid}(\mathbf{y}_n) \|_2 & \leq & 
 \left( \kappa_2(P)  \rho(A) \right)^n . \| \mathrm{wid}(\mathbf{y}_0)\|_2 
+  \frac{\left( \kappa_2(P) \rho(A) \right)^n -1}{\kappa_2(P) \rho(A) -1} . \| \mathrm{wid}(\mathbf{b}) \|_2, \\[4mm]
\end{array}
\]

This inequality puts in evidence the influence of the condition number of \(P\), the matrix of eigenvectors.
For instance, in the ideal case where the eigenvectors form an orthonormal basis, no overestimation occurs.

\subsubsection{Use of the QR Factorization}
\label{sec:QR}
When the orthogonal change of coordinates involves \(Q\) from the QR-factorization of \(A\), the algorithm can be written as
\[
\begin{array}{rcl}
A  & =  & QR, \\
x_{n+1} & = & Q y_{n+1} \\
\Leftrightarrow y_{n+1} & = & Q' x_{n+1} \\
y_{n+1} & = & Q' A Q y_n + Q' b
\end{array}
\]
In exact arithmetic, one should get
\[ y_{n+1} = RQ y_n + Q'b.\]

The interval counterpart is
\[
\begin{array}{lcl}
{\mathbf x}_{n+1} & = & Q  {\mathbf y}_{n+1} \\
 {\mathbf y}_{n+1} & = & Q' A Q  {\mathbf y}_n + Q' {\mathbf b}.
\end{array}
\]

\subsubsection{Use of the SVD Factorization}
\label{sec:SVD}
Our second and third  proposals consist in using respectively \(U\) and \(V\) from the SVD-factorization  of \(A\): from \(A = U \Sigma V'\), we use either \(B=U\) or \(B=V'\), which yields
\\
\begin{tabular}{c|c}
\begin{minipage}[c]{0.47\textwidth}
\[
\begin{array}{rcl}
x_{n+1} & = & U y_{n+1} \\
\Leftrightarrow y_{n+1} & = & U' x_{n+1}, \\
y_{n+1} & = & U' A U y_n + U' b .
\end{array}
\]
In exact arithmetic, this corresponds to
\[ y_{n+1} = \Sigma VU y_n + U'b. \]
The interval counterpart is
\[
\begin{array}{lcl}
{\mathbf x}_{n+1} & = & U  {\mathbf y}_{n+1} \\
 {\mathbf y}_{n+1} & = & U' A U  {\mathbf y}_n + U' {\mathbf b}.
\end{array}
\]
\end{minipage}
&
\begin{minipage}[c]{0.47\textwidth}
\[
\begin{array}{rcl}
x_{n+1} & = & V' y_{n+1} \\
\Leftrightarrow y_{n+1} & = & V x_{n+1}, \\
y_{n+1} & = & V A V' y_n + V b .
\end{array}
\]

In exact arithmetic, this corresponds to
\[ y_{n+1} = VU \Sigma y_n + Vb. \]
The interval counterpart is
\[
\begin{array}{lcl}
{\mathbf x}_{n+1} & = & V'  {\mathbf y}_{n+1} \\
 {\mathbf y}_{n+1} & = & V A V'  {\mathbf y}_n + V {\mathbf b}.
\end{array}
\]
\end{minipage}
\\
\end{tabular}

\noindent {\bf Remark:} \(VU\) is also an orthogonal matrix.

\subsection{Lohner's QR Method}
\label{sec:LohnerQR}

A well-known approach is given in Lohner, e.g. in \cite{Lohner2001} and studied in details by Nedialkov and Jackson in \cite{NedialkovJackson2001}.
It is usually presented for the iteration \( x_{n+1} = A_n x_n + b_n\), that is when the matrix and the affine term vary at each iteration.

Lohner's QR method consists in performing the following iteration:
\[ \left\{ \begin{array}{lcl}
\multicolumn{3}{l}{y_0=x_0, \: Q_0=I,  \: \lbrack Q_1,R_1 \rbrack =qr(A) \mbox{ that is, } A=Q_1R_1} \\
\lbrack Q_{n+1},R_{n+1} \rbrack & = & qr(R_n Q_n) \\
y_{n+1} & = & Q_{n+1}'A Q_n y_n+ Q_{n+1}'b \\
x_{n+1} & = &Q_{n+1} y_{n+1}
\end{array}
\right.
\]
and its interval counterpart is
\[ \left\{ \begin{array}{lcl}
\multicolumn{3}{l}{ {\mathbf y}_0={\mathbf x}_0, \: Q_0=I,  \: \lbrack Q_1,R_1 \rbrack =qr(A) \mbox{ that is, } A=Q_1R_1} \\
\lbrack Q_{n+1},R_{n+1} \rbrack & = & qr(R_n Q_n) \\
 {\mathbf y}_{n+1} & = & Q_{n+1}'A Q_n  {\mathbf y}_n+ Q_{n+1}'{\mathbf b} \\
{\mathbf x}_{n+1} & = &Q_{n+1}  {\mathbf y}_{n+1}. \\
\end{array}
\right.
\]

In the case of a constant -- throughout the iterations -- matrix \(A\), one can recognize Francis' and  Kublanovskaya's QR-algorithm.
Using the convergence of \( (R_n)\) towards the matrix of eigenvalues of \(A\) (or towards its Schur form),
in \cite{NedialkovJackson2001}, Nedialkov and Jackson established the following bounds:
%\[ \begin{array}{lcl}
% w({\mathbf x}_n) &  \leq & \mathrm{cond}(P) \rho(A)^n w({\mathbf x}_0) \\[2mm]
%& &  + \frac{\mathrm{cond}(P) \rho(A)^{n-1}-1}{\mathrm{cond}(P) \rho(A) -1} w({\mathbf b}) \\[2mm]
%& &  + {\mathbf b}
%\end{array}
%\]
\[ w({\mathbf x}_n) \leq \mathrm{cond}(P) \rho(A)^n w({\mathbf x}_0) + \frac{\mathrm{cond}(P) \rho(A)^{n-1}-1}{\mathrm{cond}(P) \rho(A) -1} w({\mathbf b}) + {\mathbf b} \]
where we recall \(A\) diagonalizable: \(A = P^{-1} \Lambda P \). \\

\subsection{Comparison}
Two aspects are compared: the complexity and the accuracy, that is, the bounds on the widths of the iterates, of each method.
\\

Let us first examine the computational complexity.
Let us recall that the QR-factorization, resp. SVD-factorization, of a \( d \times d\) matrix has a computational complexity of \( {\cal{O}}(d^3) \).
In the algorithms of Sections~\ref{sec:QR} and \ref{sec:SVD}, the factorization of a matrix is performed only once, and not at every iteration:
for \(n\) iterations, these algorithms thus have complexity \( {\cal{O}} (d^3 + n d^2)\).
In comparison,  Lohner's QR method has complexity
\( {\cal{O}} (n d^3)\), which is significantly larger when \(d\) is large. In comparison, the cost of the factorization is negligible when the number \(n\) of iterations is large.
\\

Let us now compare the accuracy of these different methods, from a theoretical point of view.
The bounds we get on the width of the iterate \({\mathbf x}_n\) are larger than the  bounds obtained by Nedialkov and Jackson, as the condition number of the matrix \(P\) appears to the \(n\)-th power in the formula for the QR- and SVD-algorithms, whereas it appears without this \(n\)-th power in the bound for Lohner's QR-algorithm.
As a condition number is always larger or equal to \(1\), this means that the bound for Lohner's QR-algorithm is tighter than the bounds for the QR- and SVD-algorithms.

%\[ \begin{array}{lcl}
%w( {\mathbf y}_n) &  \leq & \left( \mathrm{cond}(P) d \rho(A))\right)^n E w( {\mathbf y}_0) \\
%& &  + \frac{\left( \mathrm{cond}(P) d \rho(A) \right)^{n-1}-1}{\mathrm{cond}(P) d \rho(A) -1} \|w({\mathbf b})\| e \\
%& &  + \|w({\mathbf b})\| e
%\end{array}
%\]
%where \(A\) is diagonalizable: \(A = P^{-1} \Lambda P \), \\
%\(E\) the matrix of \(1\)s and \(e\) the vector of \(1\)s.
%Indeed, the condition number of the diagonalizing matrix \(P\) now appears to the \(n\)-th power.

\section{Experiments}

\subsection{Experimental Setup}
\label{sec:setup}
After the results on the widths of the iterates of the toy example given in Section~\ref{sec:toy-example},
Section~\ref{sec:expe-toy} presents the computation of each corner of the initial box, to illustrate that it is possible to get tight enclosures, on such a small example.

All algorithms presented in this paper, namely the naive (or brute-force) application of the iteration,
the QR-algorithm of Section~\ref{sec:QR}, the two versions of the SVD-algorithm of Section~\ref{sec:SVD},
and Lohner's QR algorithm given in Section~\ref{sec:LohnerQR} have been implemented in Octave using
Heimlich' interval package \cite{Heimlich2016},
then in Matlab using Rump's Intlab package \cite{Ru99a}.
Two other methods have been implemented and compared.
The first technique~\cite{Revol2004} consists in the determination of \(k\) such that \( \rho(|A^k|) <1\), then it computes only one iterate every \(k\) step, in other words it computes
\[ x_{(k+1)n} = A^k x_{kn} + \sum_{i=0}^{k-1} A^i b:\]
this iteration converges even when interval arithmetic is employed.
The other technique is the use of affine arithmetic, as advocated by Rump in a private communication.
In the experimental results presented below, each technique is associated to a color: \\
\begin{center}
\begin{tabular}{l|l}
{\bf algorithm} & {\bf color} \\
\hline
brute force & black \\
QR & cyan \\
SVD \(U\) & red \\
SVD \(V\) & magenta \\
Lohner's QR & dark blue \\
every \(k\)-th iterate & green \\
affine arithmetic & yellow \\
\end{tabular}
\end{center}

The factorizations use only the 
basic QR and SVD factorizations available in Matlab, but neither the pivoted QR recommended by Lohner in \cite{Lohner2001} nor more elaborate versions presented by Higham in \cite{Higham2000}.
\\

Sections~\ref{sec:expe10x10} and \ref{sec:expe100x100} contain the evolution of the radii of the iterates computed by these different techniques, for two matrix dimensions: \(10 \times 10\) and \(100 \times 100\). The \(y\)-axis for the radii uses a logarithmic scale.
For both dimensions, four kinds of matrices \(A\) have been used for the experiments. On the one hand, matrices which are well-conditioned (with a condition number of order \(10^2\)) and ill-conditioned (with a condition number of order \( 10^{10}\)) have been generated. On the other hand, the scaling of the matrices varies: matrices which are well-scaled and matrices which are ill-scaled, with the order of magnitude of their coefficients varying between \(1\) and \(10^{10}\). These are only orders of magnitudes, as the matrices, originally generated by a call to Matlab's {\tt randsvd} were then added to a multiple of the identity matrix and multiplied by a constant, in order to satisfy both \( \rho(A)<1\) and \( \rho(|A|) >1\).
It can also be noted that degrading the scaling of the matrix also degrades its condition number; in other words, a ``well-conditioned ill-scaled'' matrix has a much worse condition number than a ``well-conditioned well-scaled'' matrix, even if the required condition numbers, in the call to {\tt randsvd}, are initially the same.

All experiments have been performed on a 2.7 GHz Quad-Core Intel Core i7 with 16GB RAM.
Timings are averaged over 100 executions, except for affine arithmetic where at most 10 executions were performed.

\subsection{Toy Example}
\label{sec:expe-toy}
\begin{figure}[h]
\begin{center}
\begin{tabular}{cc}
\hspace*{-5mm} \includegraphics[width=65mm]{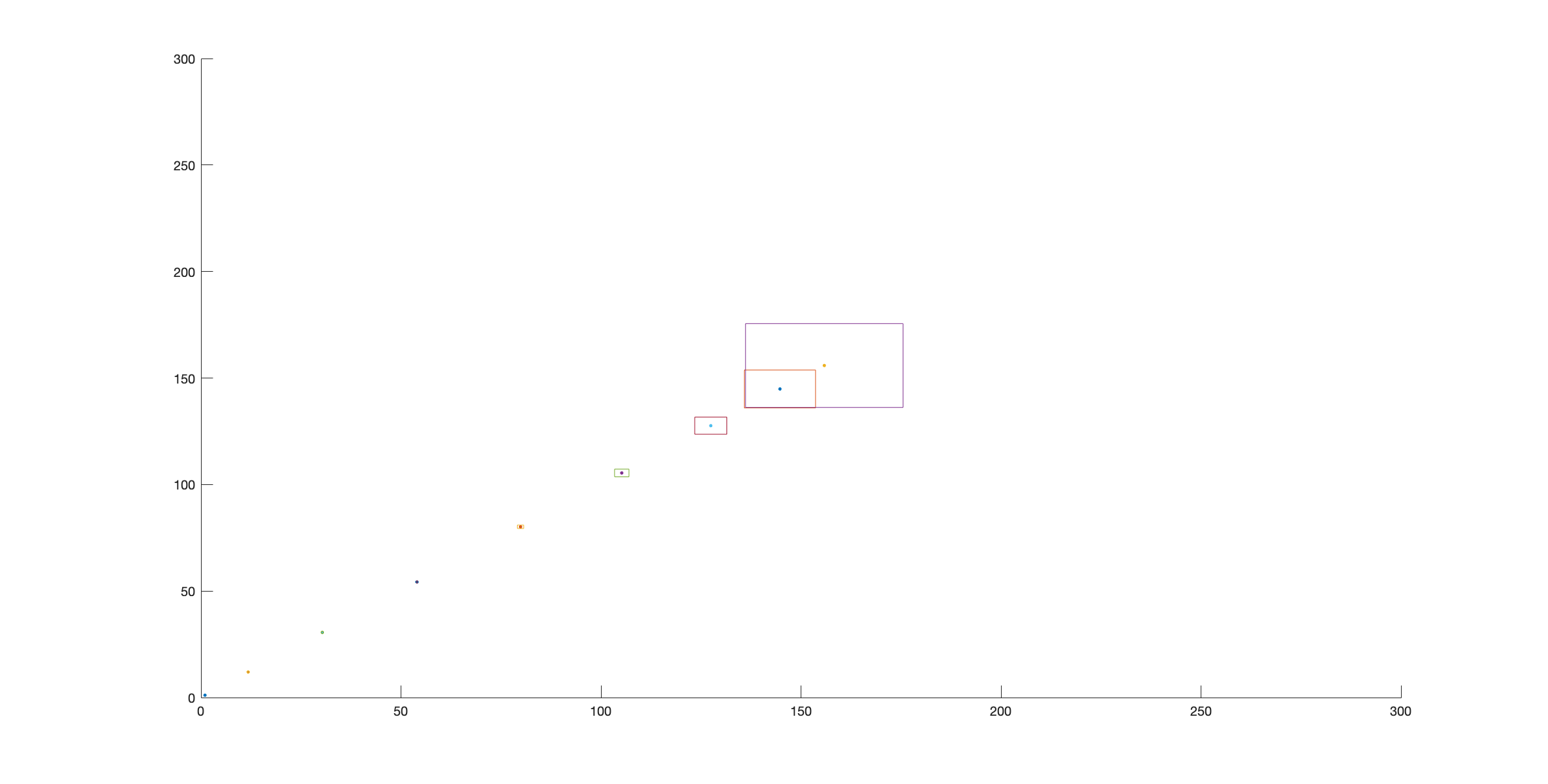}
&
\includegraphics[width=65mm]{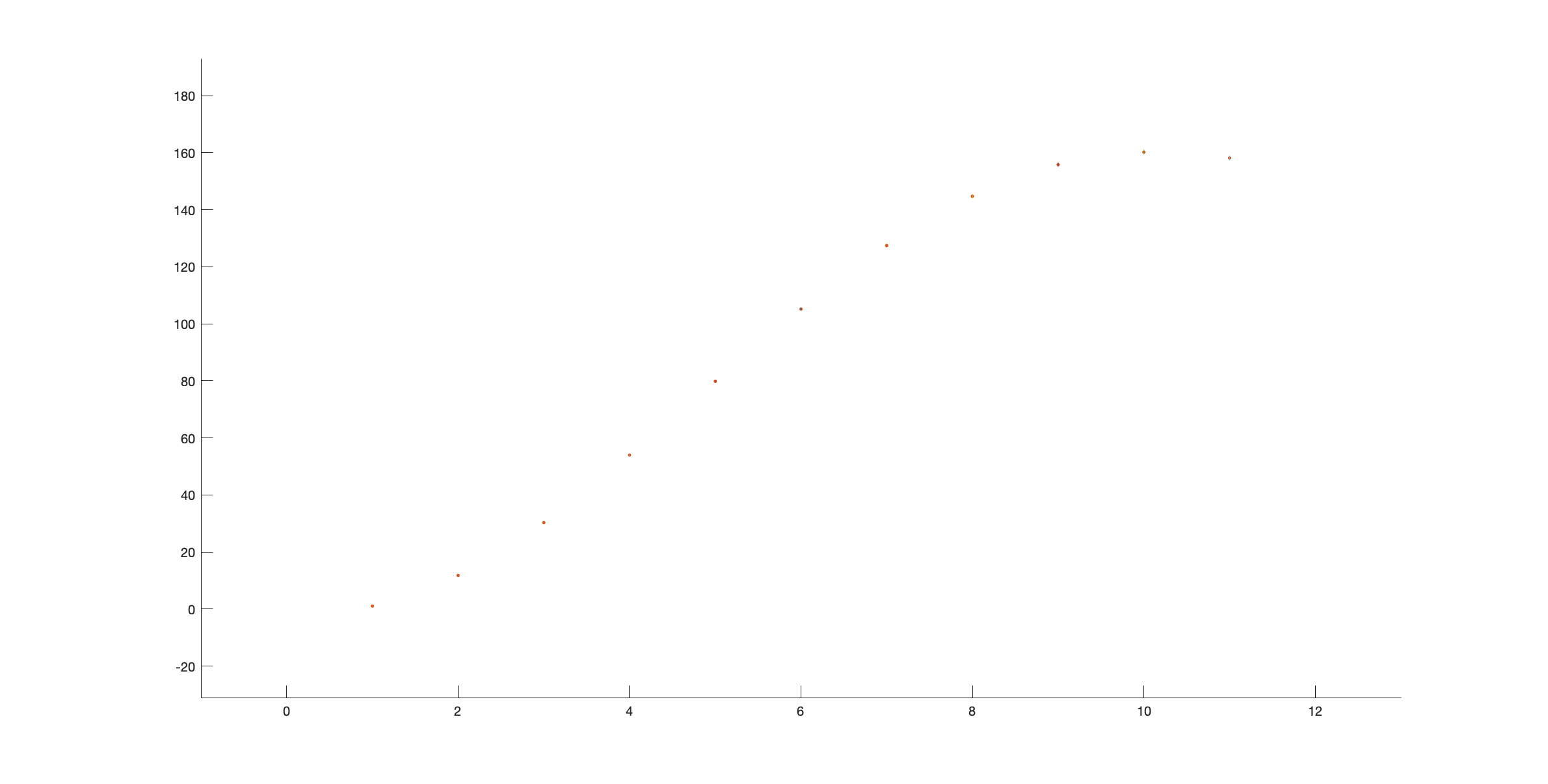}
\end{tabular}
\caption{Left: the 10 first iterates of the toy example, where the endpoints of \({\mathbf x}_0\) are considered separately.
Right: the 10 first iterates of the toy example, computed corner by corner. \label{fig:toy-example}}
\end{center}
\end{figure}

First, the toy example presented in Section \ref{sec:toy-example} is considered.
As the iteration is affine, one can compute separately  the images of the endpoints of the  initial vector, to get the endpoints of the successive iterates. That is, we compute separately
\[
x_{n} = 1.8*x_{n-1} - 0.9*x_{n-2} + 4.7.10^{-2}*3*{\mathbf u}
\]
for \( x_0 = 0\) and \(x_1 =1\) and for \( x_0 = 0\) and \(x_1 = 1.1\).
However, we use the interval vector \( {\mathbf u} = [9.95,  10.05]\) in the iteration.
The convex hull  of the 10 first iterates are represented on the left part of Figure~\ref{fig:toy-example}. It is obvious that the width of the successive iterates grow rapidly.

Then we compute separately
\[
x_{n} = 1.8*x_{n-1} - 0.9*x_{n-2} + 4.7.10^{-2}*3*u
\]
for \( x_0 = 0\) and \(x_1 =1\) and for \( x_0 = 0\) and \(x_1 = 1.1\),
and for \(u = 9.95\) and \(u = 10.05\).
The convex hull  of the 10 first iterates are represented on the right part of Figure~\ref{fig:toy-example}. In  this case, the width of the successive iterates remain small,  of the order of magnitude of \(1\%\) of the midpoint of the interval.

In this toy example, the iterations had to be performed 4 times, that is, once for each corner of the initial values, in order  to get a tight enclosure. This can  clearly not be generalized to high dimensions, as the number of corners grows as \(4^d\) with the dimension  \(d\)  of the  problem.

\subsection{Example of dimension 10}
\label{sec:expe10x10}

\begin{figure}[h]
\begin{center}
%\vspace*{-40mm}
\begin{tabular}{lcc}
a) & \includegraphics[width=65mm]{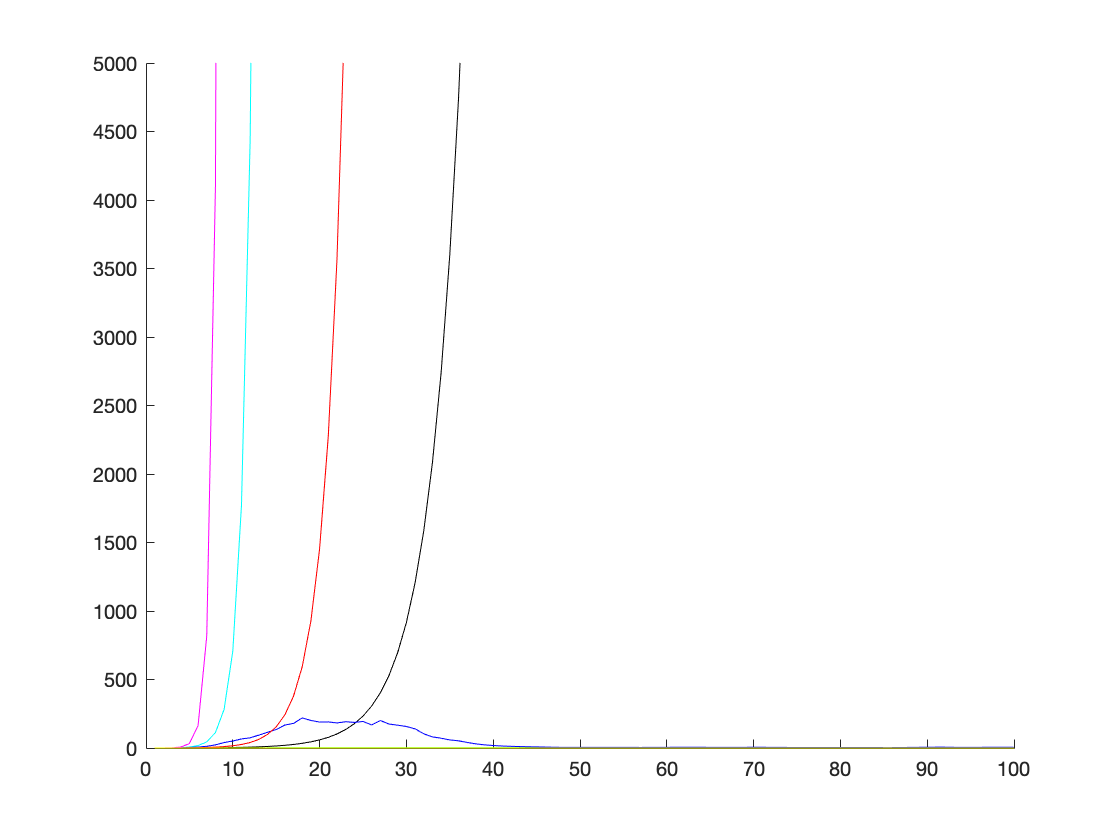}
&  
\includegraphics[width=65mm]{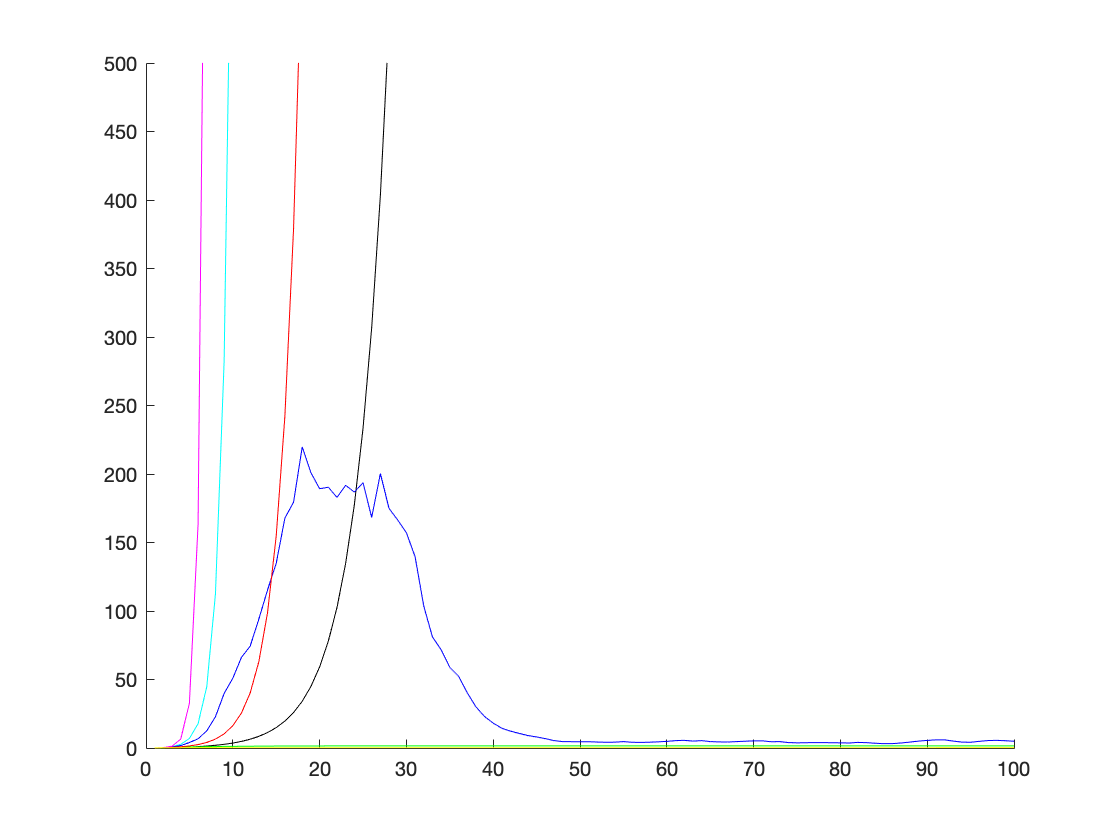}
\\
b) & \includegraphics[width=65mm]{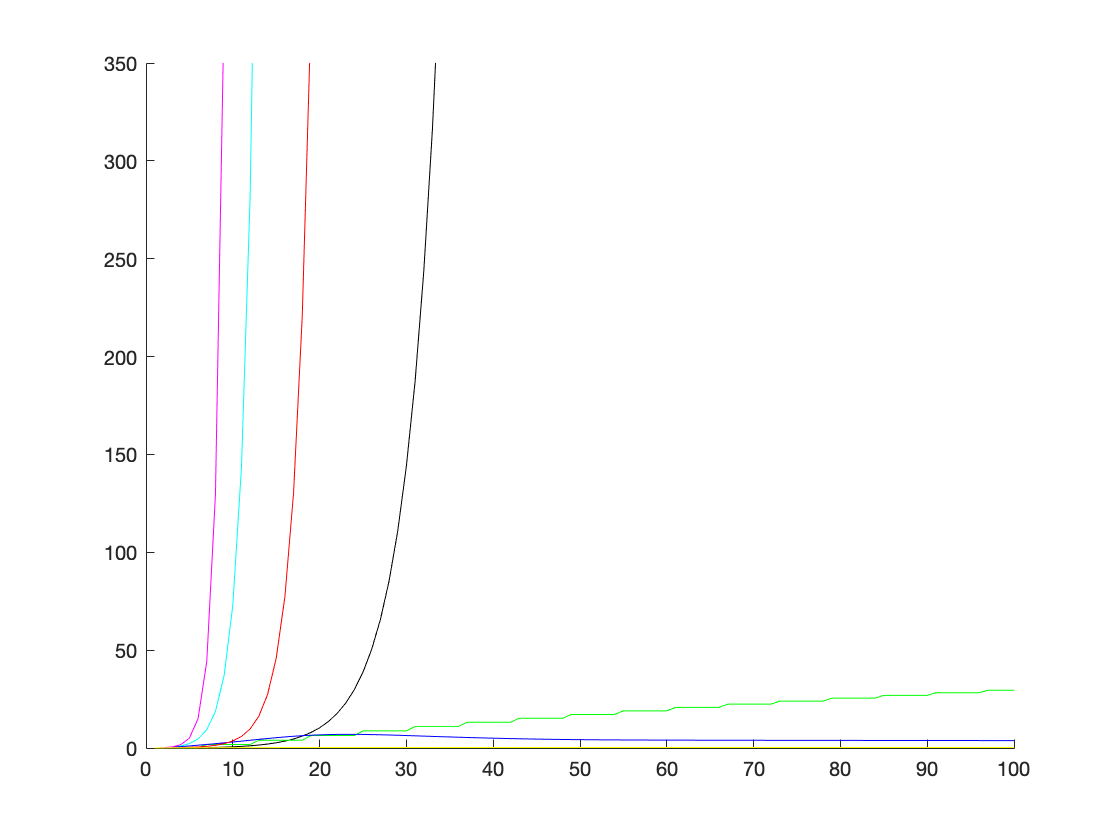}
&  
\includegraphics[width=65mm]{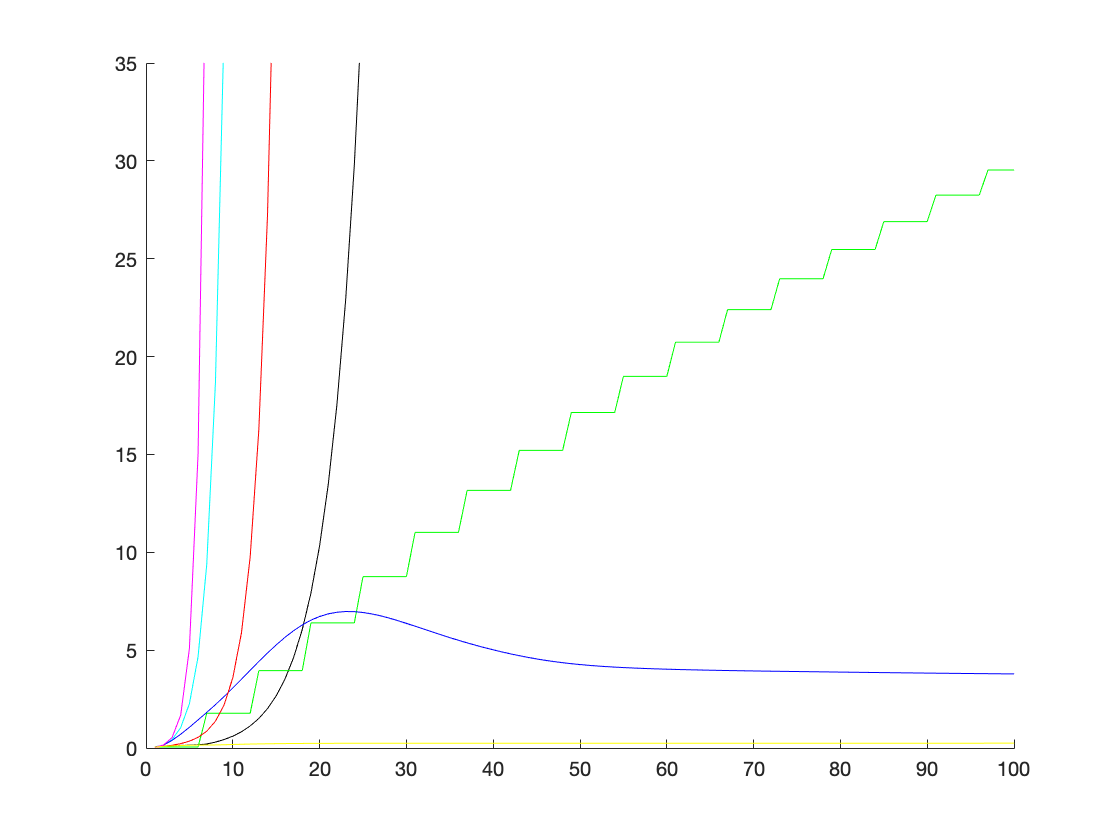}
\\
c) & \includegraphics[width=65mm]{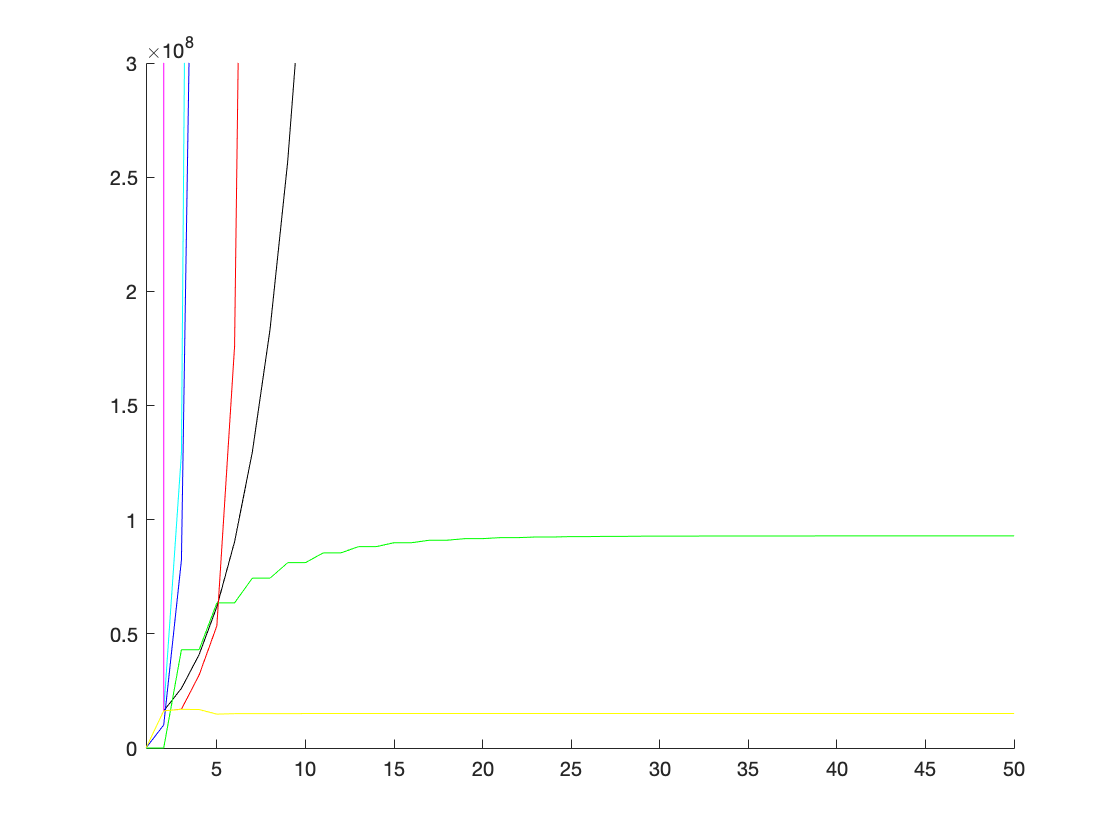}
&  
\includegraphics[width=65mm]{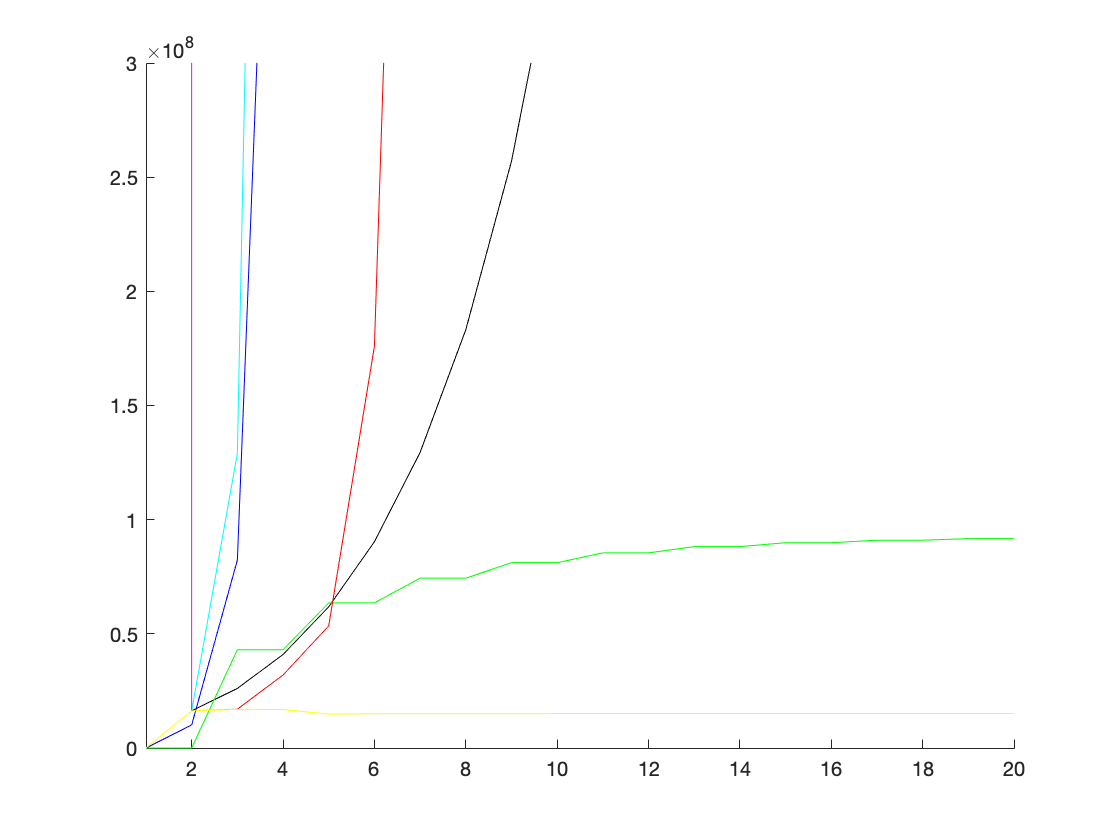}
\\
d) & \includegraphics[width=65mm]{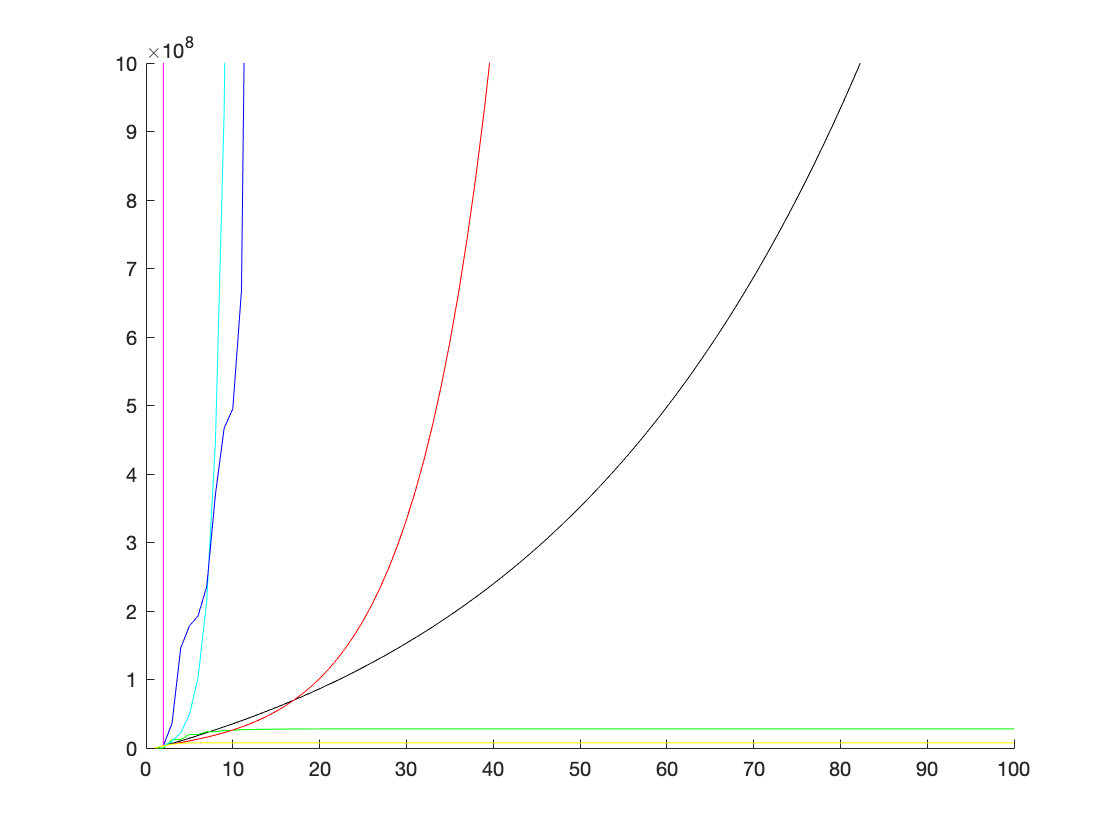}
&  
\includegraphics[width=65mm]{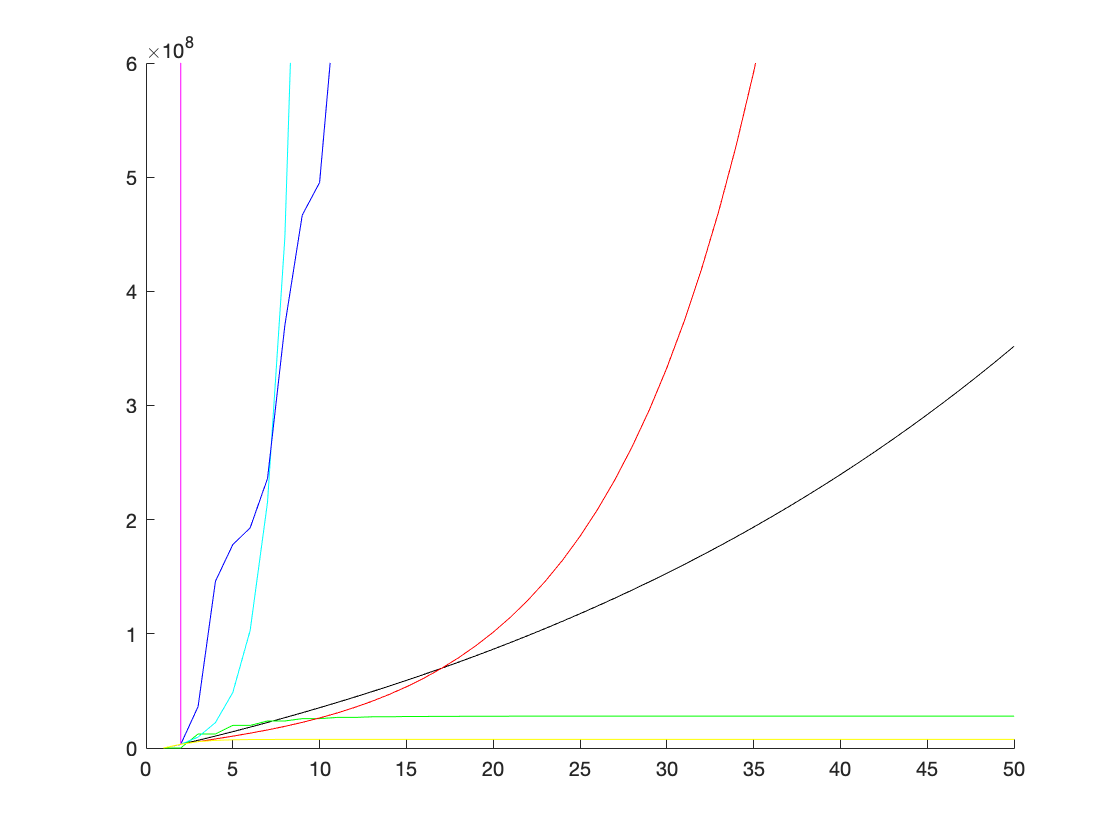}
\end{tabular}
\end{center}
\caption{The case of a \(10 \times 10\) matrix (right part: zoom of the left part): a) well-conditioned and well-scaled, b)  ill-conditioned well-scaled, c) well-conditioned ill-scaled, d) ill-conditioned ill-scaled. \label{fig:example10x10}}
\end{figure}

Figure~\ref{fig:example10x10} gives the radii (in logarithmic scale) for the successive iterates computed by the methods detailed above. 
When the number of iterations is large (visually, above 30 or 40 iterations), the iterates computed by all methods presented in Section~\ref{sec:methods} diverge rapidly, as can be seen on the plots on the left. When one concentrates on the first iterations,  the behaviours compare differently. One can also note that unscaling the matrix \(A\) speeds the divergence, for all methods.
On the contrary, the \(k\)-step method and the use of affine arithmetic preserve the convergence of the iterates.

The timings in seconds are given below:
\begin{center}
\begin{tabular}{l||r|r|r|r}
{\bf method} & {\bf well-cond.} & {\bf ill-cond.} & {\bf well-cond.} & {\bf ill-cond.} \\
& {\bf well-scaled} & {\bf well-scaled} & {\bf ill-scaled} & {\bf ill-scaled} \\
\hline
\hline
{\bf naive} & 0.0088 & 0.0084 & 0.0093 & 0.0086 \\
{\bf \(k\)-th step} & 0.0044 & 0.0015 & 0.0043 & 0.0045 \\
{\bf QR} & 0.0161 & 0.0155 & 0.0160 & 0.0157 \\
{\bf SVD \(U\)} & 0.0162 & 0.0156 & 0.0161 & 0.0166 \\
{\bf SVD \(V\)} & 0.0147 & 0.0145 & 0.0324 & 0.0332 \\
{\bf Lohner's QR} & 0.0170 & 0.0164 & 0.0165 & 0.0177 \\
{\bf affine arith.} & 8.1859 & 8.7911 & 8.6595 & 8.2754 \\
\end{tabular}
\end{center}

The methods presented in Sections~\ref{sec:QR}, \ref{sec:SVD} and \ref{sec:LohnerQR} all exhibit similar execution times.
The naive method performs less operations and is thus faster.
The \(k\)-th step method is fast as well, the variations in its execution time are due to the preprocessing, that is to the determination of the power \(k\) such that \(\rho(|A^k|) <1\): the execution time is larger when \(k\) is larger. With this method, the convergence is good.
 The use of affine arithmetic significantly slows down the computations, however the iterates converge.

\subsection{Example of dimension 100}
\label{sec:expe100x100}

\begin{figure}[h]
\begin{center}
\begin{tabular}{lcc}
a) & \includegraphics[width=65mm]{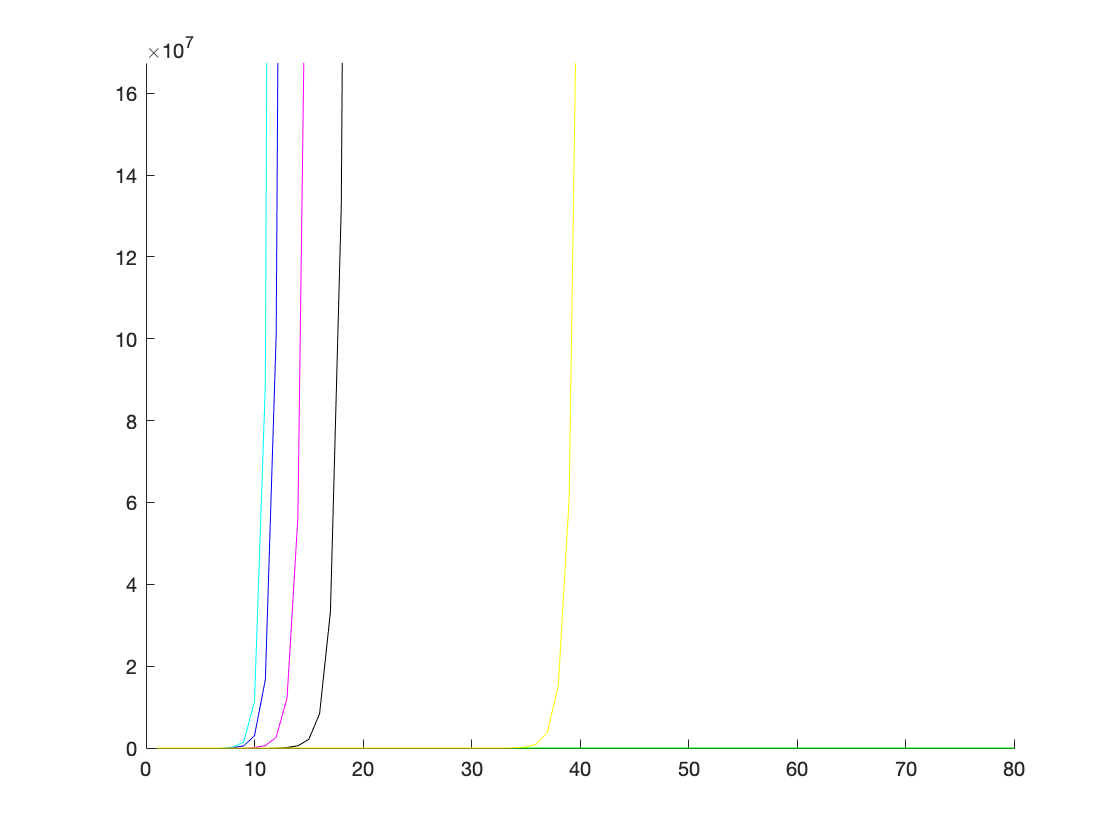}
&  
\includegraphics[width=65mm]{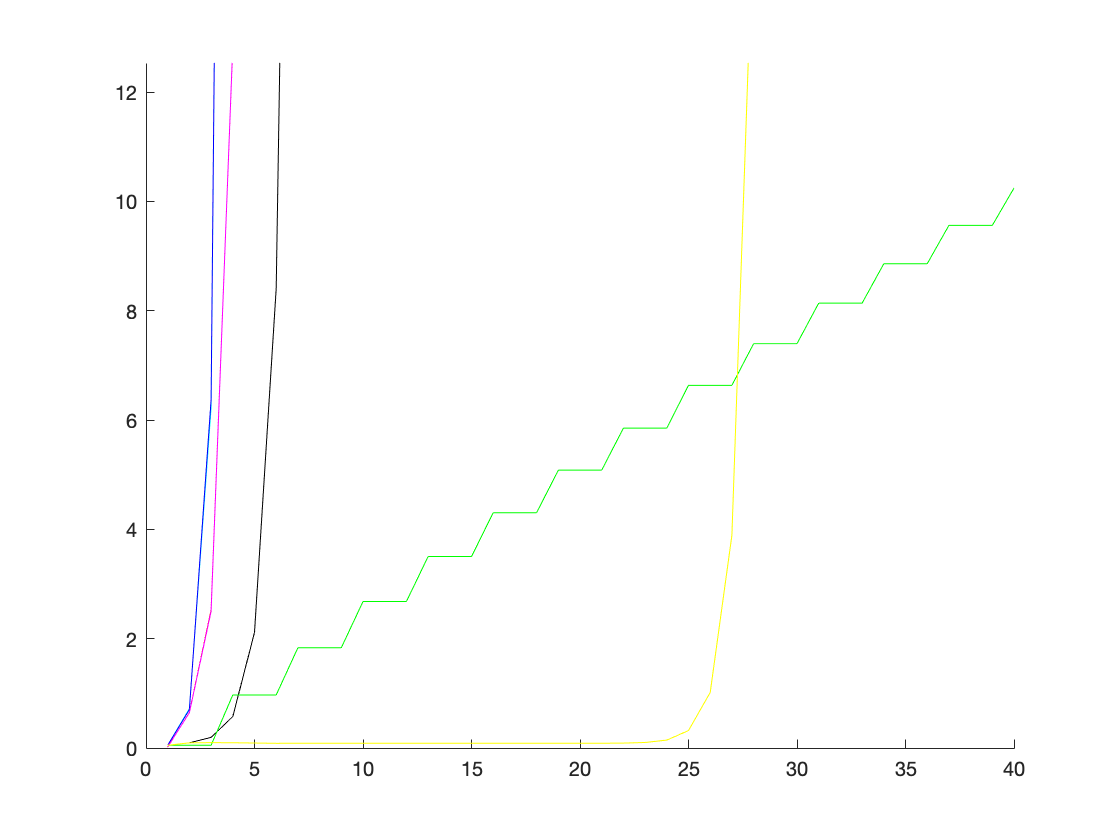}
\\
b) & \includegraphics[width=65mm]{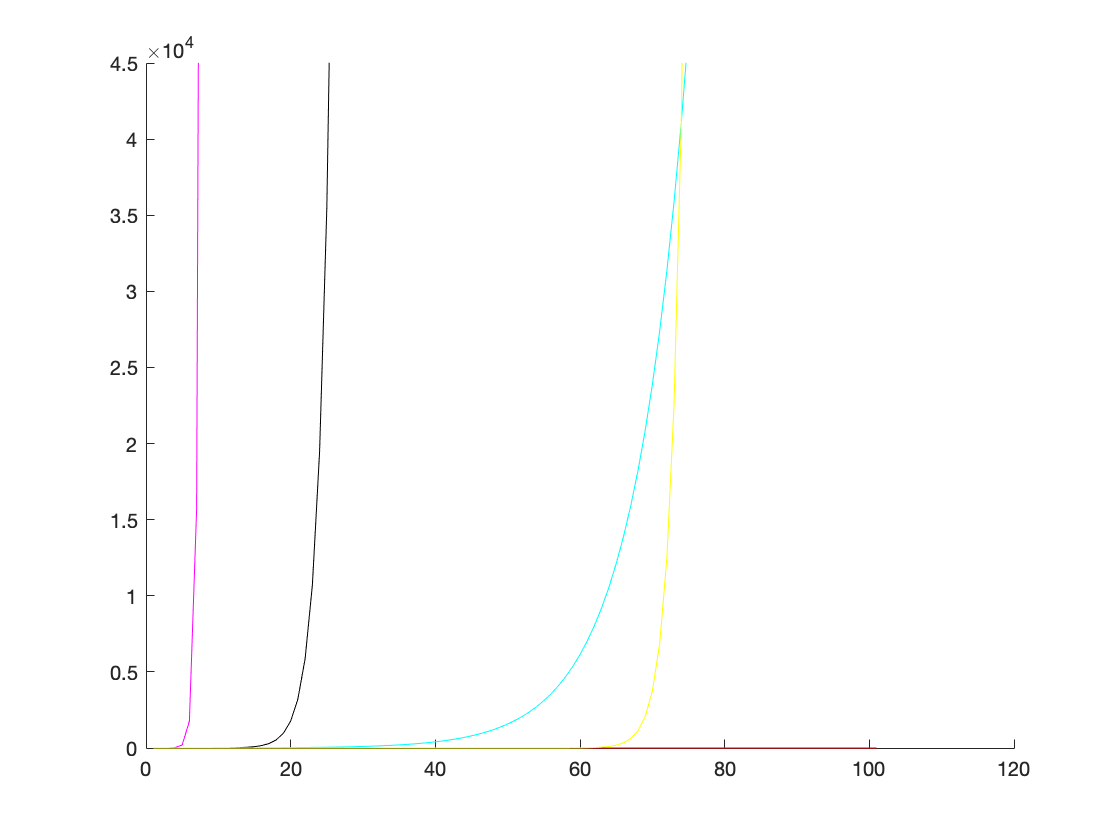}
& \includegraphics[width=65mm]{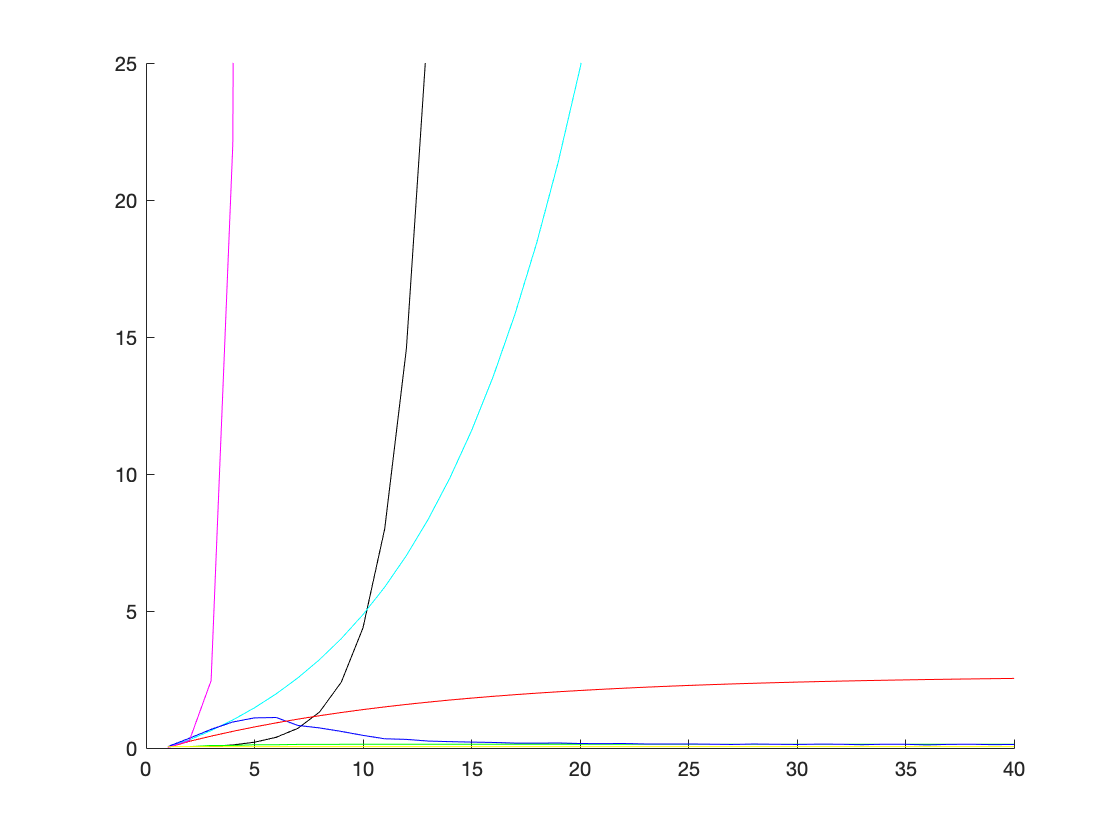}
\\
c) & \includegraphics[width=65mm]{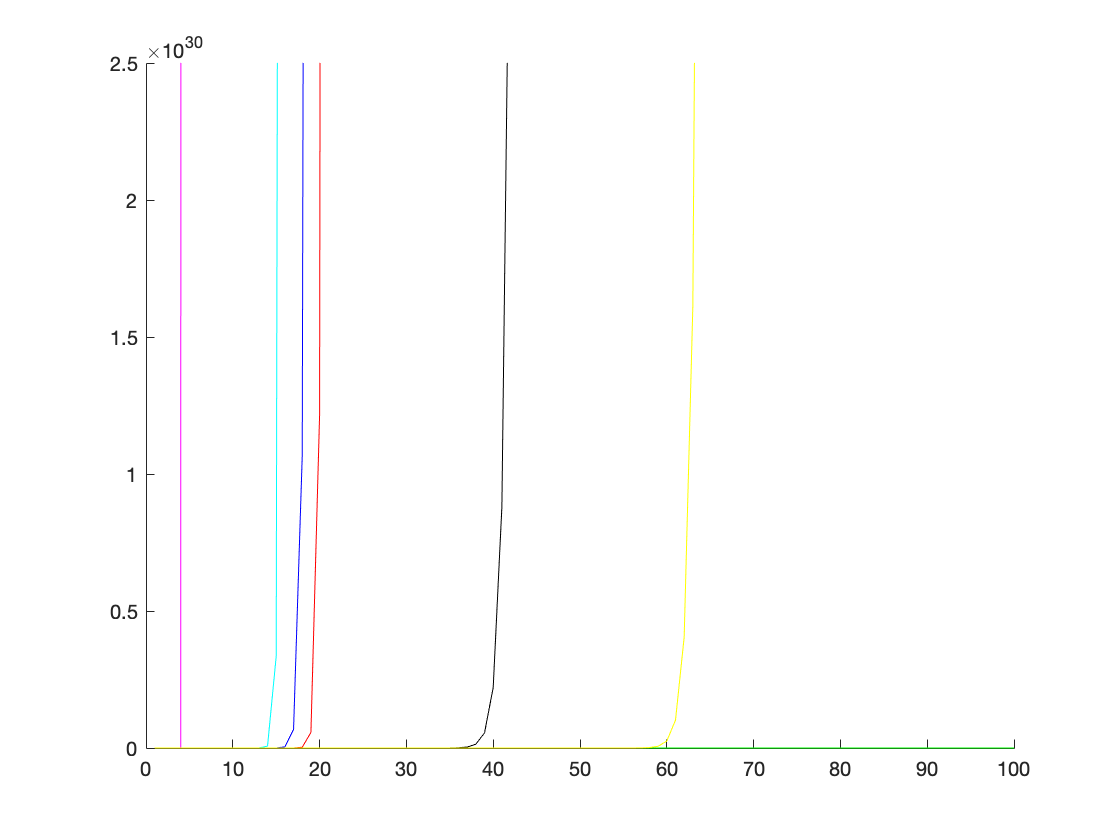}
&  
\includegraphics[width=65mm]{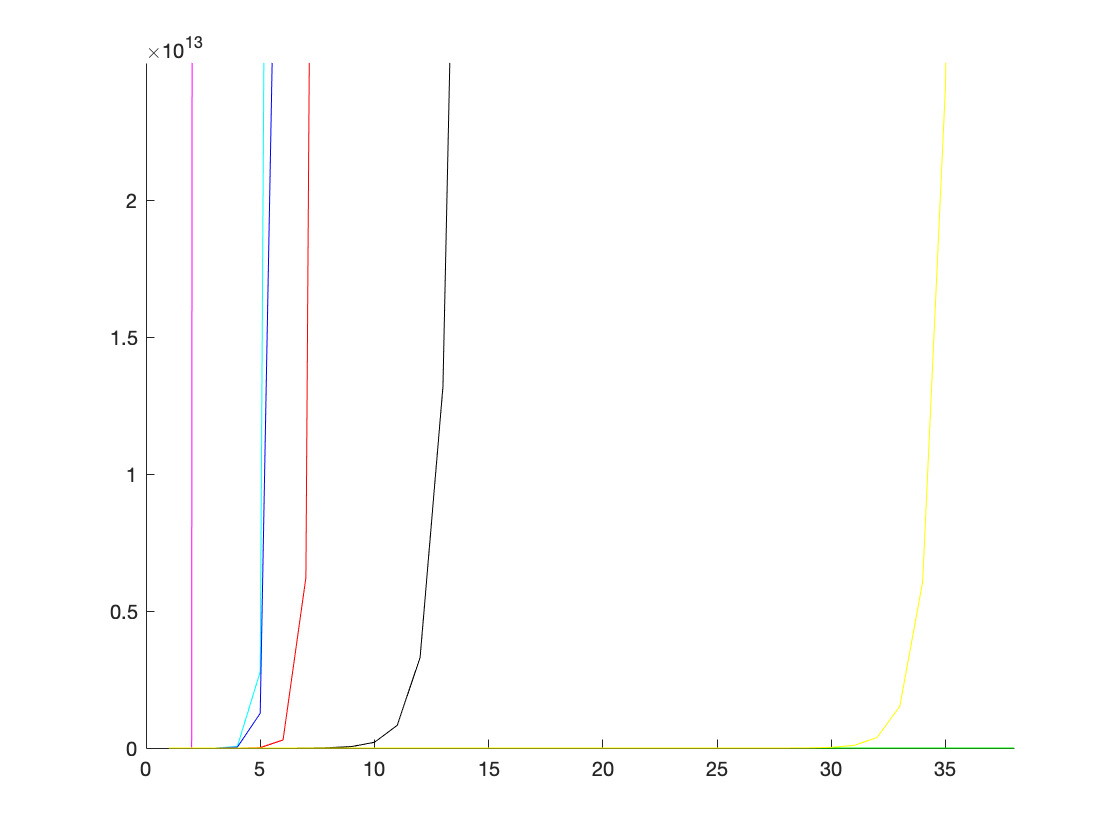}
\\
d) & \includegraphics[width=65mm]{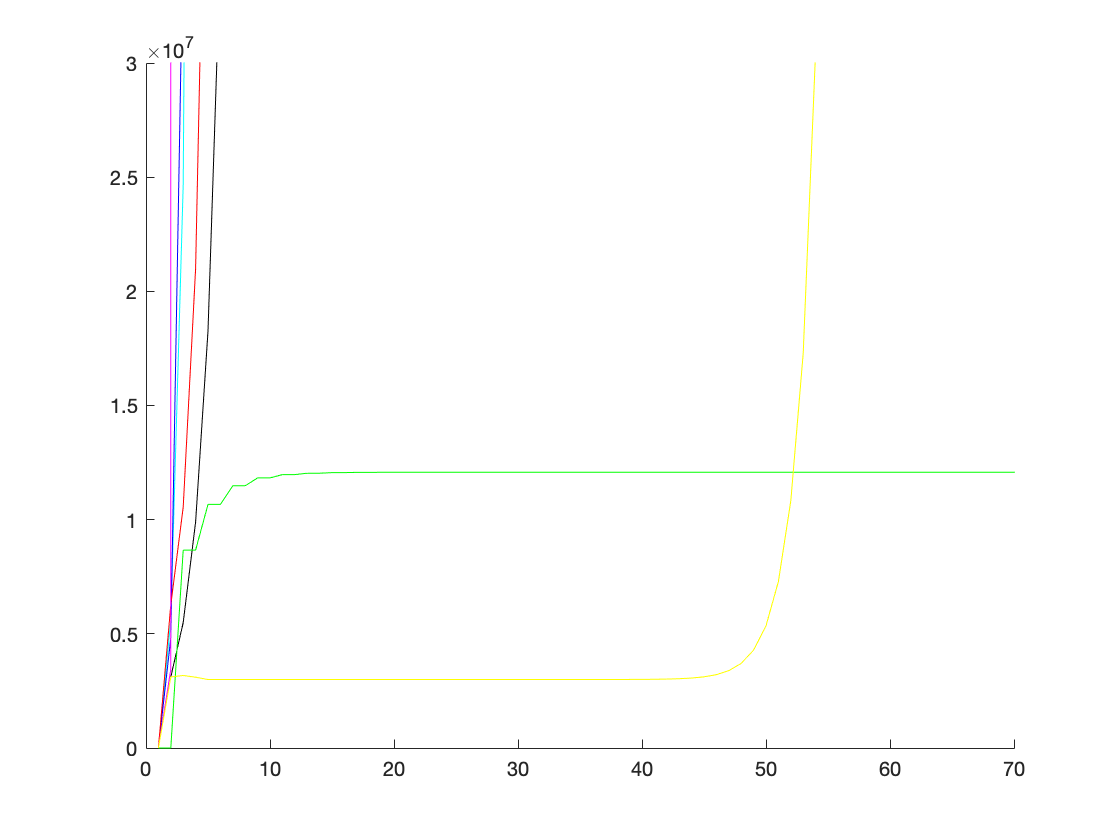}
&  
\includegraphics[width=65mm]{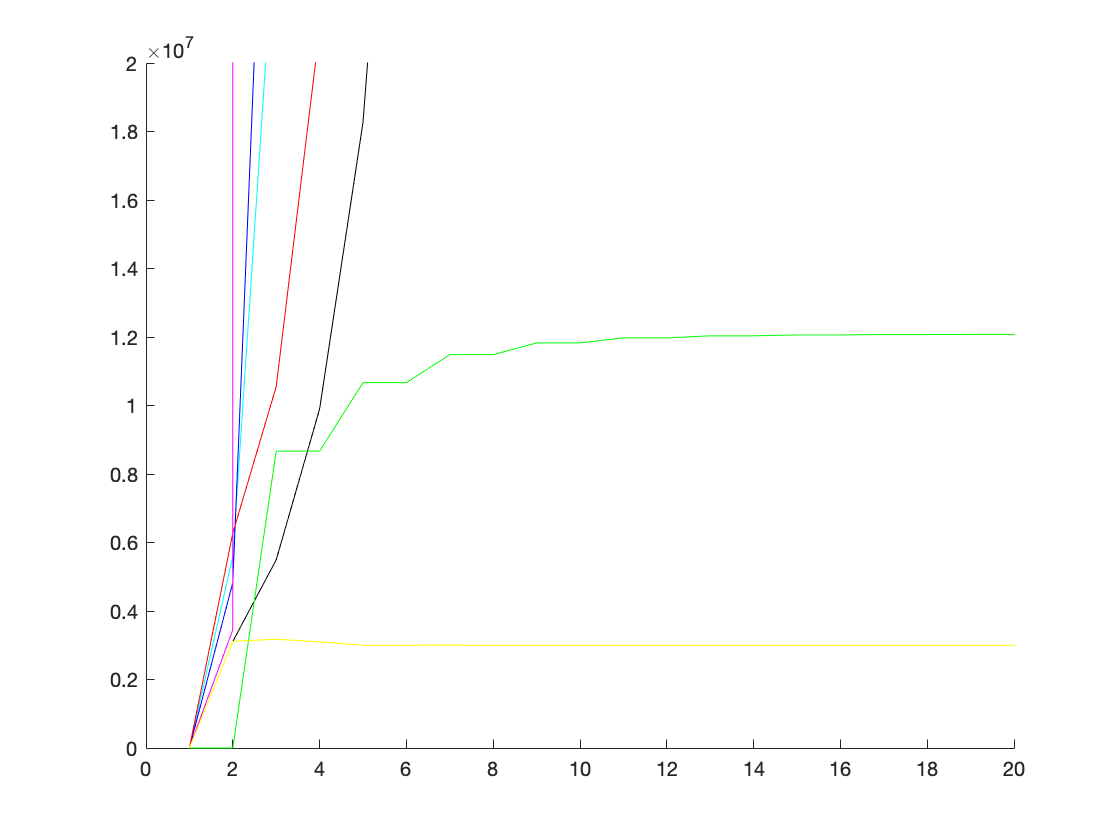}
\end{tabular}
\end{center}
\caption{The case of a \(100 \times 100\) matrix (right part: zoom of the left part): a) well-conditioned and well-scaled, b)  ill-conditioned well-scaled, c) well-conditioned ill-scaled, d) ill-conditioned ill-scaled. \label{fig:example100x100}}
\end{figure}

Figure~\ref{fig:example100x100} gives the radii (in logarithmic scale) for the successive iterates computed by the methods detailed in Section~\ref{sec:setup}. 
When the number of iterations is large (visually, above 40 or 50 iterations), the iterates computed by all methods, except the \(k\)-step method, diverge rapidly, as can be seen on the plots on the left. Again, when one concentrates on the first iterations,  the behaviours compare differently. 

The timings in seconds are given below:
\begin{center}
\begin{tabular}{l||r|r|r|r}
{\bf method} & {\bf well-cond.} & {\bf ill-cond.} & {\bf well-cond.} & {\bf ill-cond.} \\
& {\bf well-scaled} & {\bf well-scaled} & {\bf ill-scaled} & {\bf ill-scaled} \\
\hline
\hline
{\bf naive} & 0.0163 & 0.0145 & 0.0142 & 0.0142 \\
{\bf \(k\)-th step} & 0.0046 & 0.0071 & 0.0048 & 0.0072 \\
{\bf QR} & 0.1577 & 0.0363 & 0.0408 & 0.0409 \\
{\bf SVD \(U\)} & 0.0427 & 0.0420 & 0.0437 & 0.0437 \\
{\bf SVD \(V\)} & 0.0423 & 0.0390 & 0.0643 & 0.0638 \\
{\bf Lohner's QR} & 0.0611 & 0.0758 & 0.0646 & 0.0804  \\
{\bf affine arith.} & 70.8724 & 72.6441 & 76.6102 & 71.2518 \\
\end{tabular}
\end{center}

The comments on the timings apply again, with the exception of the use of affine arithmetic, which is still much slower but does not manage any more to preserve the convergence very long.

\subsection{Comments}

%TIME:
%naive: very fast,
%k-th power: most efficient both in terms of time and accuracy, however one must determine k
%
%QR and SVD : comparable in time,
%Lohner: surprisingly also comparable in time
%affine arithmetic: the "last one" to diverge, that is, small width for a much larger nb of iterations than the otheer methods, at the expense of computing time: factor much larger than 1000

%ACCURACY:
%naive: method of choice when wc-ws,
%Lohner is the method of choice when the matrix is well-scaled
%SVD (V) method of choice for n small and matrix ill-scaled, that is visible on the left part of the plot -- on the right it is superseded by other methods.

One can note that the \(k\)-step method, that is the method that resorts to a convergent interval iteration, performs very well at a moderate computation cost. Even the preprocessing time to determine the value of \(k\) has a negligible cost.

This method is a totally ad hoc approach for this problem and cannot be generalized.
However, in the framework of filters and control theory, it has a physical meaning: the divergence of the iterations  can be attributed to a sampling time which is too small to allow variations to be observed. Multiplying the sampling time by \(k\) means sampling less frequently (by a factor \(k\)) and thus being able to measure the evolution of the observed quantities.
\\

The use of affine arithmetic, on the contrary, is a very general method and it exhibits a very good accuracy, even if it eventually diverges (see the experiments with the \(100 \times 100\) matrices in Section~\ref{sec:expe100x100}).
The counterpart is the execution time, which is at least a thousand times larger than for the other methods.
This is not an issue for the experiments presented here, as the time is of order of magnitude of a minute.
\\

The methods based on the QR or SVD factorizations of the matrix \(A\)  were developed with geometric principles in mind.
For the QR-algorithm, the idea was to align the current box with the directions that are preserved by the product by \(A\), with a tradeoff between aligning the box along the eigenvectors and preserving an orthonormal system of coordinates, hence the choice of \(Q\).
For the SVD-algorithm, the idea was to align the box along the direction which gets the maximal elongation, that is along the singular vector corresponding  to  the largest singular value.

In  both cases, the benefit of these geometric transformations is mitigated  with the overestimation implied  by extra computations, and there is either  no  clear benefit for the QR-based approach, or a delicate  balance for the SVD-based approach.
The SVD-algorithm is interesting when the matrix is ill-scaled, and particularly for the first iterations.
\\

The methods of choice remain either the naive approach, when the matrix  \(A\) is well-conditioned and well-scaled, or Lohner's QR  method when the matrix is ill-conditioned. Surprisingly, the overhead of Lohner's QR method, in terms of computational time, is not as large as the formula for its complexity implies.
\\

Our general recommendation is thus:
\begin{itemize}
\item to preprocess the matrix \(A\) in order to scale it;
\item then to execute in parallel the naive approach and Lohner's QR approach,  in order to converge reasonably well for any condition number of \(A\).
\end{itemize}
Affine arithmetic is a solution of choice when other solutions fail and when the analysis and developing time is a scarce resource.

\section{Conclusion and Future Work}

%For an interesting introduction to SVD factorization and its use, see Neumaier's tutorial \cite{Neumaier1998}.

This study, both theoretical and experimental, has compared several approaches to counteract the wrapping effect for the computation of affine iterations.
Geometric considerations have led to the proposed algorithms.
The benefit of these approaches is not always clear, as a better configuration is obtained through extra-computations and thus extra-overestimation.
To deepen this geometric approach, we will aim  at simplifying the resulting formulas, at getting formulas that are closer to the mathematically equivalent, but simpler, ones that are given after each proposed transformation.
The main difficulty is to perform products such as \(Q. Q'\) or \(U'.U\), without replacing them by the identity, but in a certified  and tight way.
As the SVD-based approach seems more promising,  our future work will concentrate on the use of a certified SVD  factorization, as proposed by van der Hoeven and Yakoubsohn in~\cite{vanderHoevenYakoubsohn2018}.
We also plan to consider an interval version of the matrix, using the results in \cite{HladikDaneyTsigaridas2010} to keep guarantees on the singular quantities involved in the computations.

%\section{Bibliography}

\bibliographystyle{actaplaindoi}
\bibliography{NRevol-biblio}

\end{document}